\documentclass[11pt, twoside]{article}
\usepackage{amsmath,amsthm}
\usepackage{amsfonts}
\usepackage{amssymb,latexsym}
\usepackage{enumerate}
\newtheorem{theorem}{Theorem}[section]
\newtheorem{remark}[theorem]{Remark}
\newtheorem{proposition}[theorem]{Proposition}
\newtheorem{lemma}[theorem]{Lemma}
\newtheorem{definition}[theorem]{Definition}
\newtheorem{exam}[theorem]{Example}
\newtheorem{corollary}[theorem]{Corollary}

\def\ind{1{\hskip -3 pt}\hbox{\textsc{I}}}
\def\n{\noindent}
\def\fr{\frac}
\def\Om{\Omega}
\def\E{\mathcal E}
\def\de{\delta}
\def\ve{\varepsilon}
\def\va{\varphi}
 
\def\S{\mathcal S} 
\def\F{\mathcal F}
\DeclareMathOperator{\loc}{loc}

\begin{document}
\setlength{\baselineskip}{18truept}
\pagestyle{myheadings}
\markboth{ N. V. Phu and N.Q. Dieu}{Weighted Complex $m$-Hessian Equations}
\title {Solutions to weighted complex $m$-Hessian Equations on domains in $\mathbb C^n$}
\author{
 Nguyen Van Phu* and Nguyen Quang Dieu**
\\ *Faculty of Natural Sciences, Electric Power University,\\ Hanoi,Vietnam.\\
**Department of Mathematics, Hanoi National University of Education,\\ Hanoi, Vietnam;\\
Thang Long Institute of Mathematics and Applied Sciences, Nghiem Xuan Yem,\\ Hanoi, Vietnam
\\E-mail: phunv@epu.edu.vn\\ and ngquang.dieu@hnue.edu.vn}

\date{}
\maketitle

\renewcommand{\thefootnote}{}

\footnote{2010 \emph{Mathematics Subject Classification}: 32U05, 32W20.}

\footnote{\emph{Key words and phrases}: $m$-subharmonic functions, Complex $m$-Hessian operator, $m$-Hessian equations, $m$-polar sets, $m$-hyperconvex domain.}

\renewcommand{\thefootnote}{\arabic{footnote}}
\setcounter{footnote}{0}


\begin{abstract}
In this paper, we first study the comparison principle for the operator $H_{\chi,m}$. This result is used to solve certain weighted complex $m-$ Hessian equations.
\end{abstract}

\section{Introduction}

The complex Monge-Amp\`ere operator plays a central role in pluripotential theory and has been extensively studied through the years.  This operator was used to obtain many important results of the pluripotential theory in $\mathbb{C}^n, n>1$. In \cite{BT1} Bedford and Taylor have shown that this operator is well defined in the class of locally bounded plurisubharmonic functions with range in the class of non-negative measures. Later on, Demailly generalized the work of Bedford and Taylor for the class of locally plurisubharmonic functions with bounded values near the boundary. In \cite{Ce98} and  \cite{Ce04}, Cegrell introduced  the classes $\mathcal{F}(\Omega),\mathcal{E}(\Omega)$ which are not necessarily locally bounded and he proved that the complex Monge-Amp\`ere operator is well defined in these classes.  Recently, in \cite{Bl1} and \cite{DiKo} the authors introduced $m$-subharmonic functions which are extensions of the plurisubharmonic functions and the complex $m$-Hessian operator $H_m(.) = (dd^c.)^m\wedge \beta^{n-m}$ which is more general than the Monge-Amp\`ere operator $(dd^c.)^n$.  In \cite{Ch12},  Chinh introduced the Cegrell classes $\mathcal{F}_m(\Omega)$ and $\mathcal{E}_m(\Omega)$ which are not necessarily locally bounded and  the complex $m$-Hessian operator  is well defined in these classes. On the other hand, solving the Monge - Amp\`ere equation in the class of plurisubharmonic functions is important problem in pluripotential theory. In the classes of $m$-subharmonic functions, similar to the Monge-Amp\`ere equation,  the complex $m$-Hessian equation $H_m(u)=\mu$ also plays a similar role. This equation was first studied by Li \cite{Li04}. He solved the non-degenerate Dirichlet problem for this equation with smooth data in strongly $m$-pseudoconvex domains. One of its degenerate counterparts was studied by B{\l}ocki \cite{Bl1}, where he solved the homogeneous equation with continuous boundary data.  In \cite{Cu14}, Cuong provided a version of the subsolution theorem for the complex $m$-Hessian equation in smoothly bounded strongly $m$-pseudoconvex domains in $\mathbb{C}^n$.  Next, in \cite{Ch12} he solved complex $m$-Hessian equation in the case measures $\mu$ is dominated by $m-$ Hessian operator of a bounded $m-$ subharmonic function. In $\cite{HP17}$, the authors studied complex $m$-Hessian equation in the case when the measures $\mu$ is dominated
 by $m-$ Hessian operator of a  function in the class $\mathcal{E}_{m}(\Omega)$.
These results partially extend earlier results obtained in \cite{ahag} and \cite{ACCH09} for the plurisubharmonic case. 

In this paper, we are concerned with the existence and uniqueness of certain weighted complex $m$-Hessian equations on bounded $m-$hyperconvex domains $\Omega$ in $\mathbb C^n$. Our work is directly motivated by \cite{Cz10} where the author investigated the similar question
but for somewhat simpler operator acting on the Cegrell classes for plurisubharmonic function. 
Here by  weighted complex $m$-Hessian equations we solve an equation of the form 
$\chi (u(z),z)H_m (u)=\mu$ where $\chi$ is a certain positive measurable function defined on $(-\infty, 0) \times \Om$ and 
$\mu$ is a positive Borel measure on $\Om.$

The paper is organized as follows. Besides the introduction,  the paper has other four sections. In Section 2 we recall the definitions and results concerning the $m$-subharmonic functions which were introduced and investigated intensively in recent years by many authors (see \cite{Bl1}, \cite{SA12}). We also recall the Cegrell classes of $m$-subharmonic functions $\mathcal{F}_m(\Omega)$, $\mathcal{N}_{m}(\Omega)$ and $\mathcal{E}_m(\Omega)$ which were introduced and studied in \cite{Ch12} and $\cite{T19}$. In Section 3, we present a version of the comparison principle for the weighted $m-$ Hessian operator $H_{\chi,m}$. Finally, in Section 4, we used the obtained results  to study solutions to the  weighted $m-$ Hessian operator $H_{\chi,m}.$ For the existence of the solution, we manage to apply Schauder's fixed point theorem,
a method suggested by Cegrell in \cite{Ce84}. The problem is to create a suitable convex compact set and then appropriate continuous self maps. To make this work possible, we mention among other things, Lemma \ref{bd1} giving us a sufficient 
condition for convergence in $L^1 (\Om, \mu)$ of a weakly convergent sequence in $SH_m^{-} (\Om),$
where $\mu$ is a positive Borel measure that does not charge $m-$polar sets. We also discuss a sort of stability 
of solutions of the weighted Hessian equations. A main technical tool is Lemma \ref{caphin} about convergent in capacity of Hessian measures where we do not assume the sequence is bounded from below by a fixed element in $\F_m (\Om).$
\section{Preliminaries}

Some elements of pluripotential theory that will be used throughout the paper can be found in \cite{BT1}, \cite{Ce98}, \cite{Ce04}, \cite{Kl}, while elements of the theory of $m$-subharmonic functions and the complex $m$-Hessian operator can be found in \cite{Bl1}, \cite{SA12}.  Now we recall  the class of $m$-subharmonic functions introduced by B{\l}ocki in \cite{Bl1} and the classes $\mathcal{E}^0_m(\Omega)$, $\mathcal{F}_m(\Omega)$ which were introduced by  Chinh recently in \cite{Ch12}. Let $\Omega$ be an open subset in $\mathbb{C}^n$.  By $\beta= dd^c\|z\|^2$ we denote the canonical K\"ahler form of $\mathbb{C}^n$ with the volume element $dV_{2n}= \frac{1}{n!}\beta^n$ where $d= \partial +\overline{\partial}$ and $d^c =\frac{\partial - \overline{\partial}}{4i}$.

\noindent{\bf 2.1} First, we recall the class of $m$-subharmonic functions which were introduced and investigated in \cite{Bl1}. For $1\leq m\leq n$, we define
$$\widehat{\Gamma}_m=\{ \eta\in {\mathbb{C}}_{(1,1)}: \eta\wedge \beta^{n-1}\geq 0,\ldots, \eta^m\wedge \beta^{n-m}\geq 0\},$$
where ${\mathbb{C}_{(1,1)}}$ denotes the space of $(1,1)$-forms with constant coefficients.

\noindent\begin{definition}\label{dn1}{\rm Let $u$ be a subharmonic function on an open subset $\Omega\subset \mathbb{C}^n$. Then $u$ is said to be an {\it $m$-subharmonic} function on $\Omega$ if for every $\eta_1,\ldots, \eta_{m-1}$ in $\widehat{\Gamma}_m$ the inequality
$$ dd^c u \wedge\eta_1\wedge\cdots\wedge\eta_{m-1}\wedge\beta^{n-m}\geq 0,$$
holds in the sense of currents.}
\end{definition}

\noindent By $SH_m(\Omega)$ we denote the set of $m$-subharmonic functions on $\Omega$ while $SH^{-}_m(\Omega)$ denotes the set of negative $m$-subharmonic functions on $\Omega$. It is clear that if $u\in SH_m$ then $dd^cu\in \widehat{\Gamma}_m$.

\n Now assume that $\Omega$ is an open set in $\mathbb{C}^n$ and $u\in \mathcal{C}^2(\Omega)$. Then from the Proposition 3.1 in \cite{Bl1} (also see the Definition 1.2 in \cite{SA12}) we note that $u$ is $m$-subharmonic function on $\Omega$ if and only if $(dd^c u)^k\wedge\beta^{n-k}\geq 0,$ for $k=1,\ldots, m$. More generally, if $u_1, \ldots, u_k\in\mathcal{C}^2(\Omega),$ then for all $\eta_1, \ldots, \eta_{m-k}\in \widehat{\Gamma}_m$, we have

\begin{equation}\label{pt1}
dd^cu_1\wedge\cdots\wedge dd^cu_k\wedge\eta_1\wedge\cdots\wedge\eta_{m-k}\wedge\beta^{n-m}\geq 0
\end{equation}
holds in the sense of currents. 

\n 
We collect below basic properties of $m$-subharmonic functions that might be deduced directly from Definition 2.1. For more details, the reader may consult  $\cite{Ch15}, \cite{DH14}, \cite{SA12}.$

\begin{proposition} \label{basic}
	 Let $\Om$ be an open set in $\mathbb C^n$. Then the following assertions holds true:

\n 	
(1) If $u,v\in SH_{m}(\Omega)$ then $au+bv \in SH_{m}(\Omega) $ for any $a,b\geq 0.$
\newline
(2) $PSH(\Omega)= SH_{n}(\Omega)\subset \cdots\subset SH_{1}(\Omega)=SH(\Omega).$
\newline
(3) If $u\in SH_{m}(\Omega)$ then a standard approximation convolution $ u*\rho_{\varepsilon}$ is also an m-subharmonic function on $\Omega_{\varepsilon}=\{z\in\Omega: d(z,\partial\Omega)>\varepsilon\}$ and $u*\rho_{\varepsilon}\searrow u$  as $\varepsilon\to 0.$ 
\newline
(4) The limit of a uniformly converging or decreasing sequence of $m$-subharmonic function is  $m$-subharmonic. 
\newline
(5) Maximum of a finite number of $m$-subharmomic functions is a $m$-subharmonic function.
\end{proposition}
\noindent Now as in \cite{Bl1} and \cite{SA12} we define the complex Hessian operator for locally bounded $m$-subharmonic functions as follows.
\begin{definition}\label{dn2}{\rm Assume that $u_1,\ldots, u_p\in SH_m(\Omega)\cap L^{\infty}_{\loc}(\Omega)$. Then the {\it complex Hessian operator} $H_m(u_1,\ldots,u_p)$ is defined inductively by
$$dd^cu_p\wedge\cdots\wedge dd^cu_1\wedge\beta^{n-m}= dd^c(u_p dd^cu_{p-1}\wedge\cdots\wedge dd^cu_1\wedge\beta^{n-m}).
$$}
\end{definition}
\noindent It was shown in \cite{Bl1} and later in \cite{SA12} that $H_m(u_1,\ldots, u_p)$ is a closed positive current of bidegree $(n-m+p,n-m+p).$ Moreover, this operator is continuous under decreasing sequences of locally bounded $m$-subharmonic functions. In particular, when $u=u_1=\cdots=u_m\in SH_m(\Omega)\cap L^{\infty}_{\loc}(\Omega)$ the Borel measure $H_m(u) = (dd^c u)^m\wedge\beta^{n-m}$ is well defined and is called the complex $m$-Hessian of $u$.

\begin{exam}\label{vd1}{\rm  By using an example which is due to Sadullaev and Abullaev in \cite{SA12} we show that there exists a function which is $m$-subharmonic but not $(m+1)$-subharmonic. Let $\Omega\subset\mathbb{C}^n$ be a domain  and $0\notin\Omega$. Consider the Riesz kernel given by
$$ K_m(z) =-\frac{1}{|z|^{2(n/m-1)}}, 1\leq m < n.$$
We note that $K_m\in C^2(\Omega)$. As in \cite{SA12} we have
$$  (dd^cK_m)^k\wedge\beta^{n-k} =n(n/m-1)^k(1-k/m)|z|^{-2kn/m}\beta^n.$$
Then $(dd^cK_m)^k\wedge\beta^{n-k}\geq 0$ for all $k=1,\ldots, m$ and, hence, $K_m\in SH_m(\Omega)$. However, $(dd^cK_m)^{m+1}\wedge\beta^{n-m-1}< 0$ then $K_m\notin  SH_{m+1}(\Omega)$.}
\end{exam}

\noindent{\bf 2.2} Next, we recall the classes $\mathcal{E}^0_m(\Omega)$, $\mathcal{F}_m(\Omega)$ and $\mathcal{E}_m(\Omega)$ introduced and investigated in \cite{Ch12}. Let $\Omega$ be a bounded $m$-hyperconvex domain in $\mathbb{C}^n$, which mean there exists an $m-$ subharmonic function $\rho:\Omega\to (-\infty,0)$ such that the closure of the set $\{z\in\Omega:\rho(z)<c\}$ is compact in $\Omega$ for every $c\in (-\infty,0).$ Such a function $\rho$ is called the exhaustion function on $\Omega.$ Throughout this paper $\Omega$ will denote a bounded $m-$ hyperconver domain in $\mathbb{C}^{n}.$
 Put
$$
\mathcal{E}^0_m=\mathcal{E}^0_m(\Omega)=\{u\in SH^{-}_m(\Omega)\cap{L}^{\infty}(\Omega):
\underset{z\to\partial{\Omega}}
\lim u(z)=0, \int\limits_{\Omega}H_m(u) <\infty\},$$
$$\mathcal{F}_m=\mathcal{F}_m(\Omega)=\big\{u\in SH^{-}_m(\Omega): \exists \mathcal{E}^0_m\ni u_j\searrow u, \underset{j}\sup\int\limits_{\Omega}H_m(u_j)<\infty\big\},$$
and
\begin{align*}
&\mathcal{E}_m=\mathcal{E}_m(\Omega)=\big\{u\in SH^{-}_m(\Omega):\forall z_0\in\Omega, \exists \text{ a neighborhood } \omega\ni z_0,  \text{ and }\\
&\hskip4cm \mathcal{E}^0_m\ni u_j\searrow u \text{ on } \omega, \underset{j}\sup\int\limits_{\Omega}H_m(u_j) <\infty\big\}.
\end{align*}
\n In the case $m=n$ the classes $\mathcal{E}^0_m(\Omega)$, $\mathcal{F}_m(\Omega)$ and $\mathcal{E}_m(\Omega)$ coincide, respectively, with the classes $\mathcal{E}^0(\Omega)$, $\mathcal{F}(\Omega)$ and $\mathcal{E}(\Omega)$ introduced and investigated earlier by Cegrell in \cite{Ce98} and \cite{Ce04}.

\noindent From  Theorem 3.14 in \cite{Ch12} it follows that if $u\in \mathcal{E}_m(\Omega)$, the complex $m$-Hessian $H_m(u)= (dd^c u)^m\wedge\beta^{n-m}$ is well defined and it is a Radon measure on $\Omega$. On the other hand, by Remark 3.6 in \cite{Ch12} the following description of  $\mathcal{E}_m(\Omega)$ may be given
\begin{align*}
&\mathcal{E}_m=\mathcal{E}_m(\Omega)=\big\{u\in SH^{-}_m(\Omega):\forall U\Subset\Omega, \exists {v\in\mathcal{F}_m(\Omega)}, {v=u}  \text{ on }  {U}\bigl\}.
\end{align*}

\begin{exam}{\rm  For $0 < \alpha < 1$ we define the function
$$u_{m,\alpha} (z):= - (-\log \|z\|)^{\frac{\alpha m}{n}} + (\log 2)^{\frac{\alpha m}{n}}, 1\leq m\leq n, $$
on the ball $\Omega:= \{z\in\mathbb C^{n}:\|z\|<\frac 1 2\}$. 
Direct computations as in Example 2.3 of \cite{Ce98} shows that 
$u_{m,\alpha}\in\mathcal E_{m}(\Omega)$, $\forall 0<\alpha < \frac 1 m$. 

}\end{exam}
\noindent{\bf 2.3.} 
We say that an $m-$ subharmonic function $u$ is maximal if for every relatively compact open set $K$ on $\Omega$ and for each upper semicontinuous function  $v$ on $\overline{K},$ $v\in SH_{m}(K)$ and $v\leq u$ on $\partial K,$ we have $v\leq u$ on $K.$ The family of maximal $m-$ subharmonic function defined on $\Omega$ will be denoted by $MSH_{m}(\Omega).$ 
As in the plurisubharmonic case, if  $u\in\mathcal{E}_{m}(\Omega)$ then maximality of $u$ is characterized by $H_{m}(u)=0$ (see \cite{T19}).\\
\n {\bf 2.4.} Following \cite{Ch15}, a set $E\subset\mathbb{C}^n$ is called $m$-polar if  $E\subset \{v=-\infty\}$ for some $v\in SH_m(\mathbb{C}^n)$ and $v$ is not equivalent $-\infty.$

\noindent{\bf 2.5.} In the same fashion as the relative capacity introduced by Bedford and Taylor in \cite{BT1}, the $Cap_m$ relative capacity is defined as follows.
\begin{definition}\label{dn61}{\rm Let $E\subset \Omega$ be a Borel subset. The $m$-capacity of $E$ with respect to $\Omega$ is defined in \cite{Ch15} by
		$$Cap_m(E,\Omega)= \sup\Bigl\{\int\limits_{E}H_m(u): u\in SH_m(\Omega), -1\leq u\leq 0\Bigl\}.$$}
\end{definition}

\noindent Proposition 2.8 in \cite{Ch15} gives some elementary properties of the $m$-capacity similar to those presented in \cite{BT1}. Namely, we have:

a) $Cap_m(\bigcup\limits_{j=1}^\infty E_j)\leq \sum\limits_{j=1}^\infty Cap_m(E_j)$.

b) If $E_j\nearrow E$ then $Cap_m(E_j)\nearrow Cap_m(E)$.\\
According to Theorem 3.4 in \cite{SA12} (see also Theorem 2.24 in $\cite{Ch15}$), a Borel subset $E$ of $\Omega$ is $m$-polar if and only if $Cap_m(E)=0.$ A more qualitative result in this direction will be supplied in Corollary \ref{level}. In discussing convergence of complex Hessian operator, the following notion
stemming from the work of Xing in \cite{Xing00}, turns out to be quite useful.
\begin{definition}
	A sequence $\{u_j\} \subset SH_m (\Om)$ is said to converge in $Cap_m$ to $u \in SH_m (\Om)$ if for every $\de>0$ and every compact set $K$ of $\Om$ we have
	$$\lim_{j \to \infty} Cap_m (\{\vert u-u_j\vert>\de\} \cap K)=0.$$ 
\end{definition}
\n 
Generalizing the methods of Cegrell in \cite{Ce12}, it is proved in Theorem 3.6 of \cite{HP17} that $H_m (u_j) \to H_m (u)$ weakly if $u_j \to u$ in $Cap_m$ and if all $u_j$ are bounded from below by a fixed element of $\F_m.$ 

\n {\bf 2.6.} Let $u\in SH_{m}(\Omega),$ and let ${\Omega_{j}}$ be a fundamental sequence of $\Omega,$ which means
$\Om_j$ is strictly pseudoconvex, $\Omega_{j}\Subset\Omega_{j+1}$ and $\cup_{j=1}^{\infty}\Omega_{j}=\Omega.$ Set 
$$u^{j}(z)=\big(\sup\{\varphi(z):\varphi\in SH_{m}(\Omega), \varphi\leq u \,\,\text{on}\,\, \Omega_{j}^{c}\}\big)^{*},$$
where $\Omega_{j}^{c}$ denotes the complement of $\Omega_{j}$ on $\Omega.$\\
We can see that $u^{j}\in SH_{m}(\Omega)$ and $u^{j}=u$ on $(\overline{\Omega_{j}})^{c}.$ From definition of $u^{j}$ we see that $\{u^{j}\}$ is an increasing sequence and therefore $\lim\limits_{j\to\infty}u^{j}$ exists everywhere except on an $m-$ polar subset on $\Omega.$ Hence, the function $\tilde{u}$ defined by $\tilde{u}=\big(\lim\limits_{j\to\infty}u^{j}\big)^{*}$
is $m-$ subharmonic function on $\Omega.$ Obviously, we have $\tilde{u}\geq u.$ Moreover, if $u\in \mathcal{E}_{m}(\Omega)$ then $\tilde{u}\in  \mathcal{E}_{m}(\Omega)$ and $\tilde{u}\in MSH_{m}(\Omega).$
Set
$$\mathcal{N}_{m}=\mathcal{N}_{m}(\Omega)=\{u\in\mathcal{E}_{m}(\Omega): \tilde{u} =0.\}$$
We have the following inclusion
$$\mathcal{F}_{m}(\Omega)\subset \mathcal{N}_{m}(\Omega) \subset \mathcal{E}_{m}(\Omega).$$
Theorem 4.9 in $\cite{T19}$ shows that a function $u\in \mathcal{F}_{m}(\Omega)$ if and only if it belongs to the class $\mathcal{N}_{m}(\Omega)$ and has bounded total Hessian mass.

Let $\mathcal{K}$ be one of the classes $\mathcal{E}_{m}^{0}(\Omega), \mathcal{F}_{m}(\Omega), \mathcal{N}_{m}(\Omega), \mathcal{E}_{m}(\Omega).$ Denote by $\mathcal{K}^{a}$ the set of all function in $\mathcal{K}$ whose Hessian measures vanish on all $m-$polar set of $\Omega$. We say that a $m-$ subharmonic function defined on $\Omega$ belongs to the class $\mathcal{K}(f,\Omega)$, where $f\in \mathcal{E}_{m}\cap MSH_{m}(\Omega)$ if there exists a function $\varphi\in\mathcal{K}$ such that 
$$ f\geq u\geq f+\varphi.$$
Note that $\mathcal{K}(0,\Omega)=\mathcal{K}. $

\n 
We end this preliminary section by recalling the following H\"older type inequality proved in Proposition 3.3 of \cite{HP17}. In the case of plurisubharmonic functions, this sort of estimate was proved by Cegrell in his seminal work
\cite{Ce98}.
\begin{proposition} \label{holder}
	Let $u_{1},\cdots,u_{m}\in \mathcal{F}_{m}(\Omega).$ Then we have
	$$\int_{\Omega} H_{m}(u_{1},\cdots,u_{m})\leq \Big[\int_{\Omega} H_{m}(u_{1})\Big]^{\frac{1}{m}} \cdots \Big[\int_{\Omega}H_{m}(u_{m})\Big]^{\frac{1}{m}}.$$
\end{proposition}
\section{Comparison Principles for the Operator $H_{\chi,m}$}
Let $\chi: \mathbb{R}^{-}\times \Omega\rightarrow \mathbb{R}^{+}$ be a measurable function which is the pointwise limit of a sequence of {\it continuous} functions defined on $\mathbb{R}^{-}\times \Omega.$
The weighted $m-$Hessian operator $H_{\chi,m}$ is defined as follows
 $$H_{\chi,m}(u):= \chi(u(z),z){(dd^{c}u)}^{m}\wedge \beta^{n-m}, \ \forall u\in\mathcal{E}_{m}.$$
 Notice that this operator is well defined since $\chi (u(z), z)$ is measurable, being the pointwise limit of a sequence of measurable functions on $\Omega$.
 
The goal of this section is to presents some versions of the comparison principle for the operators $H_{m}$ and $H_{\chi,m}$. 
A basic ingredient is the following result (see Theorem 3.6 in \cite{HP17}).
Note that in the case $m=n,$ this lemma was included in Theorem 4.9 of [KH09]. We should say that all these work are rooted
in Proposition 4.2 in \cite{BT87} where an analogous result for plurisubharmonic functions may be found.
\begin{proposition} \label{hp}
Let $u, u_1, \cdots, u_{m-1} \in \E_m (\Om), v \in SH_m (\Om)$
 and $T:= ddc^u_1 \wedge \cdots dd^c u_{m-1} \wedge \beta^{n-m}.$
Then the two non-negative measures $dd^c \max (u, v) \wedge T$ and $dd^c u \wedge T$ coincide
on the set $\{v < u\}.$
\end{proposition} 
Now we start with the following versions of the comparison principle.  
\begin{lemma}\label{th3.8}
Let $u, v \in \mathcal{E}_{m}$ be such that 
\begin{equation} \label{ass}
H_m (u)=0 \ \text{on the common singular set}\ \{u=v=-\infty\}.
\end{equation}
Let $h \in SH^{-}_{m}(\Omega)$ be such that $h \geq -1$. Then the following estimate 
\begin{equation} \label{compa1}
\fr1{m!}\int\limits_{\{u<v\}} (v-u)^m (dd^ c h)^m \wedge\beta^{n-m}
\leq \int\limits_{\{u<v\}} (-h) [H_m (u)-H_m (v)]
\end{equation}
holds true if one of the following conditions are satisfies:

\n 
(a) $\liminf\limits_{z\to\partial \Omega}[u(z)-v(z)]\geq 0;$ 

\n 
(b) $u \in \F_m.$
\end{lemma}
\begin{remark} Observe that when $h=-1$ then (\ref{compa1}) reduces to the more standard form of the comparison principle
$$\int\limits_{\{u<v\}} H_m (v) \le \int\limits_{\{u<v\}} H_m (u).$$
\end{remark}
\begin{proof} We follow closely the arguments in Section 4 of \cite{KH09} where analogous results for plurisubharmonic functions are established.
	First we prove (\ref{compa1}) under the assumption (a).
By applying Lemma 5.5 in \cite{T19} to the case $k:= m, w_1=\cdots=w_k=h,$ we obtain 
\begin{align*}
&\fr1{m!}\int_{\{u<v\}} (v-u)^m (dd^ c h)^m \wedge\beta^{n-m}+\int_{\{u<v\}} (-h) (dd^c v)^m \wedge \beta^{n-m}   \\
&\leq \int\limits_{\{u<v\} \cup \{u=v=-\infty\}} (-h)(dd^c u)^m \wedge \beta^{n-m} \\
&= \int\limits_{\{u<v\}} (-h)(dd^c u)^m \wedge \beta^{n-m}
\end{align*}
Here the last line follows from the assumption (2). After rearranging these estimates we obtain (\ref{compa1}).
Now suppose (b) is true. Then for $\ve>0$ we set $v_\ve:= \max \{u, v-\ve\}.$ Then $u \le v_\ve \in \F_m$. So we may apply Lemma 5.4 in \cite{T19} to get
 \begin{equation*} 
 \fr1{m!}\int\limits_{\Om} (v_\ve-u)^m H_m (h)
 \leq \int\limits_{\Om} (-h) [H_m (u)-H_m (v_\ve)].
 \end{equation*}
which is the same as 
 \begin{equation} \label{eq9}
 \fr1{m!}\int\limits_{\{u<v-\ve\}} (v_\ve-u)^m H_m (h) 
 \leq \int\limits_{\Om} (-h) [H_m (u)-H_m (v_\ve)].
 \end{equation}
Now we apply Proposition \ref{hp} to get
$H_m (v_\ve)=H_m (u)$ on $\{u>v-\ve\}$ and $H_m (v_\ve)=H_m (v)$ on $\{u<v-\ve\}$. This yields
$$\begin{aligned} 
\int\limits_{\Om} (-h) [H_m (u)-H_m (v_\ve)]&=\int\limits_{\{u \le v-\ve\}} (-h) [H_m (u)-H_m (v_\ve)]\\
&\le\int\limits_{\{u \le v-\ve\}} (-h) H_m (u)+ \int\limits_{\{u<v-\ve\}} hH_m (v_\ve)\\
&=\int\limits_{\{u \le v-\ve\}} (-h) H_m (u)+ \int\limits_{\{u<v-\ve\}} hH_m (v).
\end{aligned}$$
Combining the above equality and (\ref{eq9}) we obtain
\begin{equation} \label{eqq3}
 \fr1{m!}\int\limits_{\{u<v-\ve\}} (v_\ve-u)^m H_m (h)+ \int\limits_{\{u<v-\ve\}} (-h)H_m (v)
 \le \int\limits_{\{u \le v-\ve\}} (-h) H_m (u).
 \end{equation} 
By Fatou's lemma we have 
$$\liminf\limits_{\ve \to 0} \int\limits_{\{u<v-\ve\}} (v_\ve-u)^m H_m (h) \ge \int\limits_{\{u<v\}} (v-u)^m H_m (h).$$
On the other hand, note that $\{u \le v-\ve\} \subset \{u<v\} \cup \{u=v=-\infty\}.$
Therefore using 
the hypothesis (\ref{ass}) we obtain
$$\lim\limits_{\ve \to 0} \int\limits_{\{u \le v-\ve\}} (-h) H_m (u)=\int\limits_{\{u<v\}} (-h)H_m (u).$$
So by letting $\ve \to 0$ in both sides of (\ref{eqq3}) we complete the proof.
\end{proof}
\n
Using the above result we are able to get useful estimates on the size of the sublevel sets of $u \in \F_m.$
\begin{corollary} \label{level}
For $u \in \F_m$ and $s>0$ we have the following estimates:
	
 \n
(i)$Cap_m (\{u<-s\}) \le \fr1{s^m} \int_{\Om} H_m (u).$

\n 
(ii) $\int\limits_{\{u \le -s\}}	H_m (u_s) \le 2^mm! \int\limits_{\{u<-s/2\}} H_m (u)$ where $u_s:=\max\{u, -s\}.$
\end{corollary}
\begin{proof}
(i)	Fix $h \in SH_m (\Om), -1 \le h<0.$ By the comparison principle Lemma \ref{th3.8} we have
	$$
	\int\limits_{\{u<-s\}} H_m (h)\le \int\limits_{\{\fr{u}s<h\}}  H_m (h)
	\le \fr1{s^m} \int\limits_{\{\fr{u}s<h\}} H_m (u)\\
	 \le \fr1{s^m} \int\limits_{\Om} H_m (u).$$
We are done.

\n 
(ii) By Lemma \ref{th3.8} we have
$$\begin{aligned} 
\int\limits_{\{u \le -s\}}	H_m (u_s) &\le \int\limits_{\{u \le -s\}} (-1-\fr{2u}s)^m	H_m (u_s)\\
&=\int\limits_{\{u \le -s\}} (-s-2u)^m H_m \Big (\max \Big \{\fr{u}s, -1 \Big\} \Big)\\
&=2^m\int\limits_{\{u \le -s\}} (-\fr{s}2-u)^m H_m \Big (\max \Big \{\fr{u}s, -1 \Big\} \Big)\\
&\le 2^m\int\limits_{\{u< -s/2\}} (-\fr{s}2-u)^m H_m \Big (\max \Big \{\fr{u}s, -1 \Big\} \Big)\\
&\le 2^m m! \int\limits_{\{u< -s/2\}} H_m(u).
\end{aligned}$$	 
The proof is thereby completed.
\end{proof}
\n 
A major consequence of Lemma \ref{th3.8} is the following version of the comparison principle 
which was essentially proved in Corollary 3.2 of \cite{ACCH09}  for the case when $m=n.$ 	
\begin{theorem}\label{th3.12}
Let $u\in \mathcal{N}_{m}(f)$ and $v\in\mathcal{E}_{m}(f).$ Assume that the following conditions hold true:

\n 
(a) $H_m (u)$ puts no mass on $\{u=v=-\infty\};$

\n 
(b) $H_{m}(u)\leq H_{m}(v)$ on $\{u<v\}.$ 

\n 
Then we have $u\geq v$ on $\Omega.$ In particular, if $H_m (u)=H_m (v)$ on $\Om$ then $u=v$ on $\Om.$
\end{theorem}
\n
Our proof below supplies more details to the original one in Corollary 3.2 of \cite{ACCH09}  for the case when $m=n.$ 
\begin{proof}
Fix $\ve>0.$ Choose $\varphi\in\mathcal{N}_{m}(\Omega)$ such that $f\geq u\geq f+\varphi$ on $\Omega.$
Let $\{\Omega_{j}\}$ be a fundamental sequence of $\Omega.$ Define
$$\varphi_{j}=\Big(\sup \{w: w\in SH_{m}(\Omega), w\leq \varphi \,\,\text{on}\,\, \Omega\setminus\overline{\Omega}_{j}\}\Big)^{*}.$$
Then $\varphi_j \in SH_m (\Omega), \varphi_j \le 0$ and $\varphi_j=\varphi$ on $\Omega\setminus\overline{\Omega}_{j}$.
This yields that 
$$\max \{u,v\} \ge v_j:=\max \{u, v+\varphi_j\} \in \mathcal E_m (\Omega).$$ 
Since $f \ge v$ on $\Omega$ we also have for every $j \ge 1$
$$\lim_{z \to \partial \Omega} (u(z)-v_j(z))=0.$$	
Now we note that (b) implies the estimate 
$$H_m (v +\varphi_j) \ge H_m (v) \ge H_m (u) \ \text{on}\ \ \{u<v\}.$$
It follows, in view of Proposition 5.2 in \cite{HP17}, that
\begin{equation} \label{eq1}
H_m (v_j) \ge H_m (u) \ \text{on}\ \ \{u<v\}.
\end{equation}	
Next, using the definition of $Cap_{m,\Omega}$ we obtain
$$\begin{aligned}
\fr{\ve^m}{m!} Cap_{m,\Om} (\{u+2\ve<v_j\})&=\fr{\ve^m}{m!} \sup \Big \{\int\limits_{\{u+2\ve<v_j\}} H_m (h): 
h \in SH_m (\Om), -1 \le h \le 0  \Big \}\\
&\le \fr1{m!}  \sup \Big \{\int\limits_{\{u+2\ve<v_j\}} (v_j-u-\ve)^mH_m (h): 
h \in SH_m (\Om), -1 \le h \le 0  \Big \}\\
&\le \fr1{m!}  \sup \Big \{\int\limits_{\{u+\ve<v_j\}} (v_j-u-\ve)^mH_m (h): 
h \in SH_m (\Om), -1 \le h \le 0  \Big \}\\
&\le  \sup \Big \{\int\limits_{\{u+\ve<v_j\}} (-h) [H_m (u)-H_m (v_j)]: 
h \in SH_m (\Om), -1 \le h \le 0  \Big \}\\
&\le \sup \Big \{\int\limits_{\{u<v\}} (-h) [H_m (u)-H_m (v_j)]: 
h \in SH_m (\Om), -1 \le h \le 0  \Big \}\\
&=0.
\end{aligned}$$
Here we apply the assumption (a) to obtain the fourth inequality and the last equality follows from (\ref{eq1})
and the inclusion $\{u+\ve<v_j\} \subset \{u<v\}.$ Thus $v_j \le u+2\ve$ outside a polar set of $\Om.$
Letting $j \to \infty$ while noting that $\varphi_j \to 0$
outside a polar set of $\Om,$ we see that $v \le  u+2\ve$ off a polar set of $\Om.$ Now subharmonicity of $u$ and $v$ forces $v \le u+2\ve$ entirely on $\Om.$ The proof is complete by letting $\ve \to 0.$
\end{proof}
\n
Using the basic properties of $m-$subharmonic functions in Proposition \ref{basic} and the comparison principle Lemma 
\ref{th3.8}, as in the plurisubharmonic case (see \cite{BT1}), we have the following quasicontinuity property of $m-$subharmonic functions (see Theorem 2.9 in \cite{Ch12} and Theorem 4.1 in \cite{SA12}).
\begin{proposition} \label{capacity}
	Let $u \in SH_m (\Om)$. Then for every $\ve>0$ we may find an open set $U$ in $\Om$ with $Cap_m (U)<\ve$
and $u|_{\Omega \setminus U}$ is continuous.	
\end{proposition}
\n
Using the above result and the Lemma \ref{th3.8}, as in the plurisubharmonic case (see \cite{BT1}), we have the following important fact about negligible sets for $m-$subharmonic functions (see Theorem 5.3 in \cite{SA12}).
\begin{proposition} \label{negli}
Let $\{u_j\}$ be a sequence of negative $m-$ subharmonic functions on $\Om.$ Set 
$u:=\sup\limits_{j \ge 1} u_j$. Then the set $\{z \in \Om: u(z)<u^* (z)\}$ is $m-$polar.	
\end{proposition} 
Now we are able to formulate a version of the comparison principle for the operator $H_{\chi,m}$ mentioned at the beginning of this section.
\begin{theorem}\label{comparison}
Suppose that the function $t \mapsto \chi (t,z)$ is decreasing in $t$ for every $z \in \Om \setminus E,$ 
where $E$ is a $m-$polar subset of $\Om.$
	Let $u\in\mathcal{N}_{m}(f), v\in\mathcal{E}_{m}(f)$ be such that $H_{\chi,m}(u)\leq H_{\chi,m}(v).$ 
	Assume also that $H_m (u)$ puts no mass on $\{u=-\infty\} \cup E.$
	Then we have $u\geq v$ on $\Om.$
\end{theorem}
\begin{proof}
We claim that $H_{m}(u)\leq H_{m}(v)$ on $\{u<v\}$. For this, fix a compact set $K \subset \{u<v\}$.
Let $\theta_j \ge 0$ be a sequence of continuous functions on $\Om$ with compact support such that 
$\theta_j \downarrow \ind_K.$ Since 
$$\chi (v,z)H_m (v) \ge \chi (u,z)H_m (u) \ \text{as measures on}\  \Om$$
we obtain
$$\begin{aligned}
\int_{\Om} \theta_j H_m(v)&=\int_{\Om} \fr{\theta_j}{\chi (v,z)} \chi(v,z)H_m (v)\\
& \ge \int_{\Om} \fr{\theta_j}{\chi (v,z)} \chi(u,z) H_m(u)\\
&=\int_{\Om} \theta_j \fr{\chi(u,z)}{\chi (v,z)} H_m (u).
\end{aligned}$$
Letting $j \to \infty$ we get 
$$\int\limits_K H_m (v) \ge \int\limits_K \fr{\chi(u,z)}{\chi (v,z)} H_m (u) \ge \int\limits_{K \setminus E} \fr{\chi(u,z)}{\chi (v,z)} H_m (u) =\int\limits_K H_m (u)$$
where the second inequality follows from the assumption that 
$\chi (u (z),z) \ge \chi (v (z),z)$ on $\{z: u(z)<v(z)\} \setminus E$ 
and the last estimate follows from the fact that $H_m (u)$ puts no mass on $E.$
Thus $H_{m}(u)\leq H_{m}(v)$ on $\{u<v\}$ as claimed.
Now we may apply Theorem $\ref{th3.12}$ to conclude $u\geq v.$	
\end{proof}
This section ends up with the following simple fact about convergence of measures where the concept of  convergence in capacity plays a role.
\begin{proposition} \label{convergence}
	Let $f, \{f_j\}_{j \ge 1}$ be quasicontinuous functions defined on $\Om$ and $\mu, \{\mu_j\}_{j \ge 1}$ be positive Borel measures on $\Om$. Then $f_j\mu_j$ converges weakly to $f\mu$ if the following conditions are satisfied:
	
	\n 
	(i) $\mu_j$ converges to $\mu$ weakly;
	
	\n 
	(ii) $f_j$ converges to $f$ in $Cap_m;$
	
	\n 
	(iii) The functions $\{f_j\}, f$ are locally uniformly bounded on $\Om;$
	
	\n 
	(iv) $\{\mu_j\}$ are uniformly absolutely continuous with respect to $Cap_m$ in the sense that 
	for every $\ve>0$ there exists $\de>0$ such that if $X$ is a Borel subset of $\Om$ and satisfies 
	$Cap_m (X)<\de$ then $\mu_j (X)<\ve$ for all $j \ge 1.$
\end{proposition}
\begin{proof}
	First we note that $\mu$ is also absolutely continuous with respect to $Cap_m$. Indeed, it suffices to apply
	(iii) and fact that for each {\it open} subset $X$ of $\Om$ we have 
	$\mu(X) \le \liminf\limits_{j \to \infty} \mu_j (X).$
	Now we let $\va$ be a continuous function with compact support on $\Om.$ Then we write
	$$\int \va [f_j d\mu_j-fd\mu] = \int \va(f_j-f)d\mu_j+ \Big [\int \va f d\mu_j-\int \va fd\mu \Big].$$
	Then using (i), (iii), (iv) and quasicontinuity of $f$ we see that the second term tends to $0$ as $j \to \infty$ while the first term also goes to $0$ in view of (ii), (iv) and (iii).
\end{proof}

\section{Weighted complex $m$-Hessian equations}
Let $\chi:  \mathbb{R}^{-}\times \Omega\to \mathbb{R}^{+} $ be a continuous function. 
Let $f\in\mathcal{E}_{m}(\Omega)\cap MSH_{m}(\Omega)$ be given. Then, under certain restriction on $\chi$ and the measure $\mu$, we have the following existence result for weighted complex $m-$Hessian equations.
\begin{theorem} \label{main1}
Let $\mu$ be a non-negative on $\Om$ with $\mu (\Om)<\infty$. Assume that the following conditions are satisfied:

\n 
(a) There exists $\va \in \mathcal F_m (f) \cap L^1 (\Om, \mu)$ 
such that $\mu \le H_m (\va);$ 

\n 
(b) $\mu$ puts no mass on $m-$polar subset of $\Om;$

\n 
(c) $\chi (t,z) \ge 1$ for all $t<0, z \in \Om.$

Then the equation
$$\chi(u,z)H_m(u)=\mu$$
has a solution $u \in \mathcal F^a_m (f) \cap L^1 (\Om, d\mu).$ Furthermore, if the function $t \mapsto \chi(t,z)$ is decreasing
for all $z$ out side a $m-$polar set then such a solution $u$ is unique.
\end{theorem}
\begin{remark}
The uniqueness of $u$ fails without further restriction on $\chi$. Indeed, consider the case $m=n$, and $\Om:=\{z: \vert z\vert<1\}.$ Let
$$u_1 (z):=\vert z\vert^2-1, u_2 (z):=\fr1{2} (\vert z\vert^2-1).$$
Set $$\Gamma_1:= \{(u_1 (z), z): z \in \Om)\}, \Gamma_2:= \{(u_2 (z), z): z \in \Om)\}.$$
 Then
$\Gamma_1 \cap \Gamma_2=\emptyset$ and $\Gamma_1 \cup \Gamma_2$ is a closed subset of $(-\infty,0) \times \Om.$
We will find a continuous function $\chi:  (-\infty,0) \times \Om \to \mathbb R$ such that
$\chi (t,z) \ge 1$ and that
\begin{equation} \label{ex}
 \chi (u_1,z)H_n(u_1)=\chi (u_2,z)H_n(u_2) \Leftrightarrow 2^n \chi (u_1(z),z)=\chi_2 (u_2 (z),z), \ z \in \Om.
\end{equation}
For this purpose, we first let $\chi=1$ on $\Gamma_1, \chi=2^n$ on
$\Gamma_2$. Next, by Tietze's extension theorem, we may extend $\chi$ to a continuous function on $(-\infty,0) \times \Omega$ such that $1 \le \chi \le 2^n$. Thus $\chi$ is a function satisfies (\ref{ex}) and of course the condition (c). 
Now we put 
$$\mu:= \chi (u_1,z)H_n(u_1)=C\chi (u_1 (z),z)dV_{2n},$$ 
where $C>0$ depends only on $n.$ So $u_1, u_2$ are two distinct solution of the Hessian equation
$\chi (u,z)H_n(u)=\mu.$
Moreover, we note that $$H_n (u_1) \le \mu \le 2^nCdV_{2n} \le H_n (C'u_1)$$ where $C'>0$ is a sufficiently large constant. 
Thus, we have shown  that $\mu$ satisfies also the conditions (a) and (b) of Theoren \ref{main1}.
\end{remark}
\n For the proof of Theorem \ref{main1}
we need the following result which is Theorem 3.7 in \cite{Gasmi}. The lemma was proved by translating the original proof in \cite{ahag} for plurisubharmonic functions to the case of  $m-$subharmonic ones.
\begin{lemma} \label{ahag}
Let $\mu$ be a non-negative, finite measure on $\Om$. Assume that $\mu$ puts no mass on $m-$polar subsets of $\Om$.
Then there exists $u \in \mathcal F_m (f)$	such that $H_m (u)=\mu$. 
\end{lemma}
The result below states Lebesgue integrable of elements in $\mathcal F_m (f).$
\begin{lemma} \label{l2}
Let $\va \in \mathcal F_m (f).$ Then $\va \in L^1 (\Om, dV_{2n}).$	
\end{lemma}
\begin{proof}
We may assume that $f=0.$	
Choose $\theta \in \E_m^0$ such that $H_m (\theta)=dV_{2n}.$ Then by integration by parts we have 
$$\int_{\Om} \va dV_{2n}= \int_{\Om} \va H_m (\theta)=\int_{\Om} \theta dd^c \va \wedge (dd^c \theta)^{m-1} \wedge \beta^{n-m}>-\infty.$$
Here the last estimate follows from H\"older inequality Proposition \ref{holder} and the fact that $\theta$ is bounded from below. 
\end{proof}
Next, we will prove a lemma which might be of independent interest.
\begin{lemma}\label{bd1}
Let $\mu$ be a positive measure on $\Omega$ which vanishes on all $m-$ polar sets and $\mu(\Omega) < \infty.$ Let $\{u_{j}\}\in SH_{m}^{-} (\Om)$ be a sequence satisfying the following conditions:

\n 
(i) $\sup\limits_{j \ge 1} \int\limits_{\Om} -u_jd\mu <\infty;$

\n 
(ii) $u_j \to u \in SH_m^{-} (\Om)$ a.e. $dV_{2n}.$

Then we have 
$$\lim_{j \to \infty} \int_{\Om} \vert u_j- u \vert d\mu=0.$$
\end{lemma}	
\n 
The above result is implicitly contained in the proof of Lemma 5.2 in \cite{Ce98}. We include the proof here only for the reader convenience. Notice that we also use some ideas in \cite{DH14} at the end of the proof of the lemma.
\begin{proof}
We split the proof into two steps.

\n 
{\it Step 1.} We will prove
\begin{equation} \label{eq3}
	\lim\limits_{j\to \infty}\int_{\Omega}u_{j}d\mu=\int_{\Omega}ud\mu.
	\end{equation}
To see this, we note that, in view of (i), by passing to a subsequence we may achieve that 
\begin{equation} \label{eq7} \lim\limits_{j\to \infty}\int_{\Omega}u_{j}d\mu=a.
\end{equation}
	Notice that, by monotone convergence theorem, we have 
$$\lim_{N \to \infty} \int_{\Om} \max \{u, -N\} d\mu=\int_{\Om} ud\mu,$$	
	and for each $N \ge 1$ fixed 
	$$\lim_{j \to \infty} \int_{\Om} \max \{u_j, -N\} d\mu=\int_{\Om} \max \{u, -N\} d\mu.$$
Therefore, using a diagonal process, it suffices to prove (\ref{eq3})	under the restriction that $u_j$ and $u$ are all uniformly bounded from below.
Since $\mu(\Om)<\infty$  we see that the set $A:=\{u_j\}_{j \ge 1}$ is bounded in
the Hilbert space $L^2 (\Om,\mu)$. Thus, by Mazur's theorem, we can find a sequence $\tilde u_j$ belonging to the convex hull of $A$ that converges to some element $\tilde u \in L^2 (\Om, \mu).$ After switching to a subsequence we may assume that $\tilde u_j \to \tilde u$ a.e. in $d\mu.$ But by (ii) $\tilde u_j \to u$ in $L^2 (\Om, dV_{2n})$ so 
$(\sup\limits_{k \ge j} \tilde u_k)^* \downarrow u \ \text{entirely on}\  \Om.$
Thus, using monotone convergence theorem we obtain
$$
\int_{\Om} ud\mu=\lim_{j \to \infty} \int_{\Om} (\sup\limits_{k \ge j} \tilde u_k)^* d\mu
=\lim_{j \to \infty}\int_{\Om} (\sup\limits_{k \ge j} \tilde u_k) d\mu
=\int_{\Om} \tilde ud\mu =a.$$
Here the second equality follows from the fact that $\mu$ does not charge the $m-$polar negligible set 
$(\sup\limits_{k \ge j} \tilde u_k)^* \ne (\sup\limits_{k \ge j} \tilde u_k),$ and the last equality results from the choice of $\tilde u_j$ and (\ref{eq7}). The equation (\ref{eq3}) follows.

\n 
{\it Step 2.} Completion of the proof. Set $v_j:= (\sup\limits_{k \ge j} u_k)^*$. Then $v_j \ge u_j, v_j \downarrow u$ on $\Omega$ and $v_j \to u$ in $L^1 (\Om, dV_{2n}).$
So by the result obtained in Step 1 we have
\begin{equation} \label{eq4}
\lim_{j \to \infty} \int_{\Om} v_j d\mu= \int_\Om ud\mu=\lim_{j \to \infty} \int_{\Om} u_j d\mu. 
\end{equation} 
Using the triangle in equality we obtain
$$\begin{aligned} 
\int_{\Om} \vert u_j-u\vert d\mu &\le \int_{\Om} (v_j-u)d\mu+ \int_{\Om} (v_j-u_j)d\mu\\
&=  2\int_{\Om} (v_j-u)d\mu + \int_{\Om} (u-u_j)d\mu.
\end{aligned}$$
Hence by applying (\ref{eq4}) we finish the proof of the lemma.
\end{proof}
\n
Now, we  turn to the proof of Theorem \ref{main1} where the fixed point method from \cite{Ce84} will be crucial.
\begin{proof} (of Theorem \ref{main1})
We set
$$\mathcal A:=\{u \in \mathcal F_m (f): \va \le u \le f \}.$$
First using Lemma \ref{l2} we see that $\mathcal A$ is a compact convex subset of $L^1 (\Om, dV_{2n}).$
Moreover, from the assumption on $\mu,$ and Lemma \ref{bd1} we infer that $\mathcal A$ is also compact in $L^1 (\Om, \mu).$
Let $\mathcal S: \mathcal A \to \mathcal A$ be the operator assigning each element $u \in \mathcal A$ to the {\it unique} solution $v:=\mathcal S(u) \in \mathcal F_m (f)$ of the equation
$$H_m (v)= \fr1{\chi (u(z),z)} d\mu.$$ 
This is possible according to Lemma $\ref{ahag}$, because by (b), the measure on the right hand side does not charge $m-$polar subsets of $\Om$. Note also that for such a solution $v \in \mathcal F_m (f),$ by (a) and (c), we have $H_m (v) \le \mu \le H_m (\va)$.
So the comparison principle (Theorem \ref{th3.12}) yields that $v \ge \va$ on $\Om.$ Hence the operator $\mathcal S$ indeed maps $\mathcal A$ into itself.
The key step is to check continuity (in $L^1 (\Om)$) of $\mathcal S$. Thus, given a sequence $\{u_j\}_{j \ge } \subset \mathcal A, u_j \to u$ in $L^1 (\Om)$. We must show $\S(u_j) \to \S(u)$  in $L^1 (\Om)$.
By passing to subsequences of $u_j$ coupling with Lemma $\ref{bd1}$,  we may assume that $u_j \to u$ a.e. ($d\mu$). Now we define for $z \in \Om$ the following sequences of non-negative bounded measurable functions
$$\psi^1_j (z):= \inf_{k \ge j} \fr1{\chi (u_k (z),z)},
\psi^2_j (z):= \sup_{k \ge j} \fr1{\chi (u_k (z),z)}.$$
Then we have:\\
\n 
(i) $0 \le \psi^1_j (z) \le \fr1{\chi (u_j (z),z)} \le \psi^2_j (z) \le 1$ for $j \ge 1;$

\n 
(ii) $\lim\limits_{j \to \infty} \psi^1_j (z)=\lim\limits_{j \to \infty} \psi^1_2 (z)=\fr1{\chi (u(z),z)}$ a.e. ($d\mu$).

Now, using Lemma \ref{ahag} we may find $v^1_j, v^2_j \in \mathcal F_m (f)$ are solutions of the equations 
$$H_m (v^1_j)=\psi^1_j d\mu, H_m (v^2_j)=\psi^2_j d\mu.$$
Then, using  the comparison principle we see that $v^1_j \downarrow v^1, v^2_j \uparrow v^2$, furthermore, in view of (i) we also have
\begin{equation} \label{eq2}
v^1_j \ge S(u_j) \ge v^2_j.
\end{equation}
Next we use (ii) to get 
$$H_m (v^1_j) \to \fr1{\chi (u,z)} d\mu, H_m (v^2_j) \to \fr1{\chi (u,z)} d\mu.$$
So by the monotone convergence theorem we infer 
$$H_m (v^1)=H_m ((v^2)^*)=\fr1{\chi{(u(z),z)}}d\mu=H_m(\mathcal S(u)).$$
Applying again the comparison principle we obtain $v^1=(v^2)^*=\mathcal S(u)$ on $\Om.$ 
By the squeezing property (\ref{eq2}), $S(u_j) \to S(u)$ pointwise outside a $m-$polar set of $\Om.$
Since $\mu$ puts no mass on $m-$polar sets, we may apply Lebesgue dominated convergence theorem
to achieve that $\S(u_j) \to \S(u)$  in $L^1 (\Om, d\mu)$. Thus $\mathcal S: \mathcal A \to \mathcal A$ is continuous. 
So we can invoke Schauder's fixed point theorem to attain $u \in \mathcal A$ such that 
$u=\mathcal S(u).$ Note also that $H_m (u)$, being dominated by $\mu,$ does not charge $m-$polar sets, so 
$u \in \F^a_m (f).$
Hence $u$ is a solution of the weighted $m-$Hessian equation that we are looking for.
Finally, under the restriction that $\chi(t,z)$ is decreasing for all $z$ out side a $m-$polar set, we may apply Theorem \ref{comparison} to achieve the uniqueness of such a solution $u.$
\end{proof}
\n
In our next result, we deal with the situation when $\mu$ is dominated by a suitable function of $Cap_m.$
This type of result is somewhat motivated from seminal work of Kolodjiez in \cite{Klo}.
\begin{theorem} \label{mainmain}
Let $\mu$ be a non-negative Borel measure on $\Om$ with $\mu(\Om)<\infty$ and 
$F: [0,\infty) \to [0,\infty)$ be non-decreasing function with $F(0)=0$ and 
\begin{equation} \label{eq10}
\int_1^\infty F(\fr1{s^m}) ds<\infty.\end{equation}
Assume that the following conditions are satisfied:

\n 
(a) $\mu (X) \le F(Cap_m (X))$ for all Borel subsets $X$ of $\Om;$

\n 
(b) There exists a measurable function $G: \Om \to [0, \infty]$ such
$$\chi (t,z) \ge G(z), \ \forall (t,z) \in (-\infty, 0) \times \Omega \ \text{and}\  
c:=\int\limits_{\Om} \fr1{G} d\mu  <\infty.$$

Then the equation
$$\chi(u,z)H_m(u)=\mu$$
has a solution $u \in \mathcal F_m \cap L^1 (\Om,\mu).$ 
\end{theorem}	
\begin{remark}
According to Proposition 2.1 in \cite{DiKo}, for every $p \in (0, \fr{n}{n-m})$ there exists a constant $A$ depending only on $p$ such that 
$$V_{2n} (X) \le ACap_m (X)^p$$
for all Borel subsets $X$ of $\Om.$ So the Lebesgue measure $dV_{2n}$ satisfies the assumption (a) for $F(x)=Ax^p$
and $p$ is any number in the interval $(\fr1{m}, \fr{n}{n-m}).$
\end{remark}
\begin{proof} 
Let 
$$\mathcal A:=\Big \{u \in \F_m: \int_{\Om} H_m (u) \le c\Big \}.$$
First, using H\"older inequality Proposition \ref{holder}, we will show $A$ is convex. Indeed,
let $\alpha\in [0,1],$ 
it suffices to prove $\int\limits_{\Omega}H_{m}( \alpha u+ (1-\alpha)v)\leq c.$
For this, we use Proposition \ref{holder} to get
\begin{align*}
\int_{\Omega}H_{m}( \alpha u+ (1-\alpha)v)&= \int_{\Omega}dd^{c}( \alpha u+ (1-\alpha)v)^{m}\wedge \beta^{n-m}\\
&=\int_{\Omega}\sum_{k=0}^{m} \binom{m}k\alpha^{k} (1-\alpha)^{m-k}(dd^{c}u)^{k} \wedge (  dd^{c} v)^{m-k}\wedge \beta^{n-m}\\
&= \sum_{k=0}^{m}\binom{m}k\alpha^{k} (1-\alpha)^{m-k}\int_{\Omega}(dd^{c}u)^{k}\wedge (dd^{c} v)^{m-k}\wedge \beta^{n-m}\\
&\leq \sum_{k=0}^{m}\binom{m}k\alpha^{k} (1-\alpha)^{m-k} \Big[\int_{\Omega}H_{m}(u)\Big]^{\frac{k}{m}}\Big[\int_{\Omega}H_{m}(v)\Big]^{\frac{m-k}{m}}\\
&\leq \Big[\sum_{k=0}^{m}\binom{m}k \alpha^{k} (1-\alpha)^{m-k}\Big] c=c.
\end{align*}
Thus we have proved that $\mathcal{A}$ is indeed convex.
We want to show $\mathcal A$ is compact in $L^1 (\Om, \mu).$ Indeed, first by Lemma \ref{l2} we have $\mathcal A \subset L^1 (\Om, dV_{2n}).$ 
Next we let $\{u_j\}$ be a sequence in $\mathcal A.$ 
By Lemma \ref{level}, for $s>0$ we have
\begin{equation} \label{eqq1}
Cap_m (\{u_j<-s\}) \le \fr1{s^m} \int_{\Om} H_m (u_j) \le \fr{c}{s^m}.
\end{equation}
So, in particular $u_j$ cannot contain converge to $-\infty$ uniformly on compact sets of $\Om$. Hence
by passing to a subsequence we may achieve that $u_j$ converges in  $L_{loc}^1 (\Om, dV_{2n})$
to $u \in SH_m (\Om), u<0.$ 
Notice that, using the comparison principle as in Lemma 2.1 in \cite{Cz10} we conclude that $u \in \mathcal F_m.$
Now we claim that $u_j \to u$ in $L^1 (\Om,\mu).$
In view of Lemma \ref{bd1}, it suffices to check that 
\begin{equation} \label{eqq2}
\sup_{j \ge 1} \int_{\Om} (-u_j)d\mu<\infty. 
\end{equation}
For this purpose, we apply (\ref{eqq1}) and the assumption (a) to obtain
$$\mu (\{u_j<-s\}) \le F(Cap_m (\{u_j<-s\})) \le F(\fr{c}{s^m}).$$
Hence
$$\sup_{j \ge 1}\int_{\Om} (-u_j)d\mu=\sup_{j \ge 1}\int_0^\infty \mu (\{u_j<-s\}) ds<\infty$$
where the last integral converges in view of (\ref{eq10}). Thus the claim (\ref{eqq2}) follows.
By Lemma \ref{bd1} we have $u_j \to u$ in $L^1 (\Om, d\mu).$ From now on, our argument will be close to that of the proof of Theorem \ref{main1}. More precisely, let $\mathcal S: \mathcal A \to \mathcal A$ be the operator assigning each element $u \in \mathcal A$ to the {\it unique} solution $v:=\mathcal S(u) \in \mathcal F_m$ of the equation
$$H_m (v)= \fr1{\chi (u(z),z)} d\mu.$$ 
This is possible according to Lemma $\ref{ahag}$, because by (a) and (b), the measure on the right hand side does not charge $m-$polar subsets of $\Om$ and has total finite mass $ \le c$. 
By repeating the same reasoning as in the proof of Theorem \ref{main1} (the only notable change is to replace the upper bound of the sequence $\{\psi^2_j\}$ by $\fr1{G}$) we can see that $\mathcal A$ is continuous. Thus, applying again Schauder's fixed point theorem we conclude that $\mathcal S$ admits a fixed point which is a solution of the equation $\chi (u,z)H_m(u)=\mu.$ 
The proof is then complete.
\end{proof}	
Our article ends up with the following "weak" stability result.
\begin{theorem} \label{stability}
Let $\Om, \mu, F, \chi$ and $G$ be as in Theorem \ref{mainmain}. Let $\mu_j$ be a sequence of positive Borel measures on $\Om$ such that $\mu_j \le \mu$ and $\mu_j$ converges weakly to $\mu.$ Let $u_j \in \F_m$ be a solution of the equation
$$\chi (u(z), u) H_m (u)=\mu_j.$$
Assume that $F$ and $\chi$ satisfies the following additional properties:

\n 
(i) $\int\limits_1^{\infty} F(\fr1{s^{2m}})ds<\infty;$

\n 
(ii) $\fr1{G} \in L^2 (\Om, d\mu);$

\n
(iii) $\mu':=\fr1{G}\mu$ is absolutely continuous with respect to $Cap_m;$

\n 
(iv) For every compact subsets $K$ of $\Om$ and $t_0 \in (-\infty, 0)$ we have:

\n 
(a) $\sup \{\chi (t,z): t<t_0, z \in K\}<\infty;$

\n 
(b) There exists a constant $C>0$ 
(depending on $K,t_0$) such that
for $t<t'<t_0$ and $z \in K$ the estimate below holds true
$$\vert \chi (t,z)-\chi (t',z) \vert \le C\vert t-t'\vert.$$
\n 
(v) $\chi$ is continuous on $(-\infty,0) \times \Om.$

\n 
Then there exists a subsequence of $u_j$ converging in $Cap_m$ to $u \in \F_m$ such that
$$\chi (u(z), u) H_m (u)=\mu.$$	
\end{theorem}
\n 
We require the following convergence result for the operator $H_m$. This is inspired from Theorem 1 in \cite{Xing00}.
\begin{lemma} \label{caphin}
	Let $\{u_j\}$ be a sequence in $\F_m$ that converges to $u \in \F_m$ in $Cap_m$. Assume that 
	\begin{equation} \label{capdep}
	\lim_{a \to \infty} \Big (\limsup\limits_{j \to \infty}  \int\limits_{\{u_j<-a\}} H_m (u_j) \Big)=0.
	\end{equation}
	Then $H_m (u_j)$ converges weakly to $H_m (u).$	
\end{lemma}
\begin{proof} 
	Fix a continuous function $\va$ with compact support in $\Om.$ For $a>0$ we set
	$$u_{j,a}:=\max\{u_j, -a\}, u_{a}:=\max\{u, -a\}.$$
	Then we have 
	$$\begin{aligned}
	\int_{\Om } \va [H_m (u_j)-H_m (u)] &=\int_{\Om} \va [H_m (u_j)-H_m (u_{j,a})]\\
	&+\int_{\Om} \va [H_m (u_{j,a})-H_m (u_{a})]\\
	&+\int_{\Om} \va [H_m (u_{a})-H_m (u)].
	\end{aligned}$$
	Note that, by Theorem 3.6 in \cite{HP17} we have $\int\limits_{\Om} \va [H_m (u_{a})-H_m (u)] \to 0$ as $a \to \infty$ and 
	$\int\limits_{\Om} \va [H_m (u_{j,a})-H_m (u_{a})] \to 0$ as $j \to \infty$ for any {\it fixed} $a>0.$ Thus it suffices to check
	\begin{equation} \label{capxinh}
	\lim_{a \to \infty} \Big (\limsup\limits_{j \to \infty} \Big \vert \int_{\Om} \va [H_m (u_j)-H_m (u_{j,a})] \Big \vert \Big)=0.
	\end{equation} 
	For this, we observe that $H_m (u_{j,a})=H_m (u_j)$ on the set $\{u_j>-a\}$ by Proposition \ref{hp}.
	 It now follows, using Corollary \ref{level} (ii), that
	$$\begin{aligned} 
	\Big \vert \int_{\Om} \va [H_m (u_j)-H_m (u_{j,a})] \Big \vert&= \Big \vert \int\limits_{\{u_j \le -a\}} \va [H_m (u_j)-H_m (u_{j,a})] \Big \vert\\
	& \le \Vert \va \Vert_{\Om} \Big [\int\limits_{\{u_j \le -a\}} H_m (u_j)+ \int\limits_{\{u_j \le -a\}} H_m (u_{j,a})\Big]\\
	&\le (2^m m!+1) \Vert \va \Vert_{\Om} \int\limits_{\{u_j<-a/2\}} H_m (u_j).
	\end{aligned}$$
	Thus (\ref{capxinh}) follows immediately from the assumption (\ref{capdep}). We are done.
\end{proof} 
\begin{proof}
Since $$\int\limits_{\Om} H_m (u_j) \le \int\limits_{\Om} \fr1{G} d\mu_j \le \int\limits_{\Om} \fr1{G} d\mu<\infty, \ \forall j$$
by Lemma \ref{l2}, the sequence $\{u_j\}$ is bounded in $L^1 (\Om, dV_{2n}).$
Thus after switching to a subsequence we may assume $u_j$ converges in $L^1 (\Om, dV_{2n})$ to 
$u \in SH_m (\Om).$
Our main step is to check that $u_j \to u$ in $Cap_m.$
To this end, set $\mu':=\fr1{G} \mu$, we will first claim that $u_j \to u$ in $L^1 (\Om,\mu').$
Since $\mu$ and hence $\mu'$ puts no mass on $m-$polar sets, in view of Lemma \ref{bd1}, it suffices to show
\begin{equation} \label{cap}
\sup_{j \ge 1} \int\limits_{\Om} (-u_j)d\mu'<\infty. 
\end{equation}
For this purpose, we apply Corollary \ref{level} (i) to get
\begin{equation} \label{eqq1}
Cap_m (\{\vert u_j \vert^2>s\})=Cap_m (\{u_j<-s^{1/2}\}) \le \fr1{s^{2m}} \int_{\Om} H_m (u_j) \le \fr{\mu (\Om)}{cs^{2m}}.
\end{equation}
\n
So by the assumption (a) and (\ref{eqq1}) we obtain
$$\mu (\{\vert u_j \vert^2>s\}) \le F(Cap_m (\{u_j<-s^{1/2}\})) \le F(\fr{\mu (\Om)}{cs^{2m}}).$$
This implies
$$\sup\limits_{j \ge 1} \int_{\Om} \vert u_j \vert^2d\mu=\sup\limits_{j \ge 1} \int_0^\infty \mu (\{\vert u_j \vert^2 >s\}) ds<\infty$$
where the last integral converges in view of the assumption (i).
Hence, using Cauchy-Schwarz's inequality and the assumption (ii) we obtain (\ref{cap}).
Now we turn to the convergence in $Cap_m$ of $u_j$. Fix a compact set $K$ of $\Om$ and $\de>0$. Then  
by Lemma \ref{th3.8}, for $h \in SH_m (\Om), -1 \le h<0$, we have
$$\begin{aligned}
\int\limits_{\{u-u_j>\de\}} H_m (h) &\le (\fr2{\de})^m \int\limits_{\{u-u_j>\de\}} (u-u_j-\fr{\de}2)^m H_m (h)\\
&\le (\fr2{\de})^m \int\limits_{\{u>u_j+\fr{\de}2\}} (u-u_j-\fr{\de}2)^m H_m (h)\\
&\le  (\fr2{\de})^m \int\limits_{\{u-\fr{\de}2>u_j\}} (-h) H_m (u_j)\\
&\le  (\fr2{\de})^m \int\limits_{\{u-\fr{\de}2>u_j\}} \fr1{\chi (u_j (z),z)} d\mu_j\\
&\le  (\fr2{\de})^m \int\limits_{\{u-\fr{\de}2>u_j\}} \fr1{G}d\mu\\
&\le (\fr2{\de})^{m+1} \int\limits_{\Om} \vert u_j-u \vert d\mu'.
\end{aligned}$$
It follows that
$$Cap_m (\{u-u_j>\de\}) \le (\fr2{\de})^{m+1} \int\limits_{\Om} \vert u_j-u \vert d\mu' \to 0 \ \text{as}\ j \to \infty.$$
Here the last assertion follows from Lemma \ref{bd1}.
Thus 
$$\lim_{j \to \infty} Cap_m (\{u-u_j>\de\})=0.$$
Given $\ve>0$, by quasi-continuity of $u$ we can find an open subset $U$ of $\Om$ with $Cap_m (U)<\ve$ 
such that $u$ is continuous on the compact set $K \setminus U.$ Then by Dini's theorem for all $j$ large enough the set $\{u_j-u>\de\} \cap K$ is contained in $U.$ 
So we have
$\lim\limits_{j \to \infty} Cap_m (\{u_j-u>\de\} \cap K)=0.$
Putting all these facts together we obtain 
$$\lim\limits_{j \to \infty} Cap_m (\{|u_j-u|>\de\} \cap K)=0.$$ 
So, $u_j$ indeed converges to $u$ in $Cap_m$ as claimed. We now wish to apply Lemma \ref{caphin}.
For this, fix $a>0.$ Then we have 
$$\begin{aligned} 
\int\limits_{\{u_j<-a\}} H_m (u_j)&= \int\limits_{\{u_j<-a\}} \fr1{\chi (u_j (z),z)} d\mu_j\\
& \le \int\limits_{\{u_j<-a\}} \fr1{G} d\mu_j=\int\limits_{\{u_j<-a\}} d\mu'.
\end{aligned}$$
In view of (iii) and (\ref{eqq1}) we infer that the last term goes to $0$ uniformly in $j$ as $a \to \infty$.
Thus we may apply Lemma \ref{caphin} to reach that $H_m (u_j)$ converges weakly to $H_m (u)$.	
To finish off, it remains to check $\chi (u_j (z),z) \to \chi(u(z),z)$ in $Cap_m.$  
To see this, we use the extra assumption (iv)(b) and the fact we have proved above that 
$u_j \to u$ in $Cap_m.$
Now we are in a position to apply  Proposition \ref{convergence}. In details, we note the following facts:

\n 
(a) $\chi (u_j (z), z)$ and $\chi (u (z),z)$ are quasicontinuous on $\Om,$ since $u_j$ and $u$ are such functions and since $\chi$ is continuous on $(-\infty, 0) \times \Om$ by the assumption $(v)$;

\n 
(b) $\chi (u_j (z), z)$ and $\chi (u (z),z)$ are locally uniformly bounded on $\Om.$ To see this, it suffices to note that on each compact subset $K$ of $\Om$ the functions $\{u_j\}$ and $u$ are bounded from above by a fixed constant $t_0<0$,
so by the assumption (iv)(a) we obtained the required local uniform boundedness;

\n 
(c) The sequence $\{H_m (u_j)\},$ being dominated by $\mu',$ are uniformly absolutely continuous with respect to $Cap_m$ in view of the assumption $(iii).$

It follows that 
$$\mu_j=\chi (u_j (z),z) H_m (u_j) \to \chi (u(z),z) H_m (u)$$
weakly in $\Om$. Therefore $\chi (u(z),z) H_m (u)=\mu.$
The proof is then complete.	
\end{proof}

\end{document}